\documentclass[a4paper,11pt]{amsart}
\usepackage{amsmath,amsthm,amssymb}
\usepackage[mathscr]{eucal}
\usepackage{mathrsfs}
\usepackage{enumitem}
\usepackage[bookmarks,bookmarksnumbered]{hyperref}

\numberwithin{equation}{section}

\theoremstyle{plain}
\newtheorem{main}{Theorem}
\newtheorem{mcor}[main]{Corollary}

\newtheorem{theorem}{Theorem}[section]
\newtheorem{lemma}[theorem]{Lemma}
\newtheorem{proposition}[theorem]{Proposition}
\newtheorem{corollary}[theorem]{Corollary}
\theoremstyle{definition}
\newtheorem{definition}[theorem]{Definition}
\newtheorem{notation}[theorem]{Notation}

\newcommand{\N}{\mathbb N}
\newcommand{\R}{\mathbb R}
\newcommand{\C}{\mathbb C}
\newcommand{\F}{\mathbb F}
\newcommand{\cH}{\mathcal H}
\newcommand{\cK}{\mathcal K}
\newcommand{\cL}{\mathcal L}
\newcommand{\cU}{\mathcal U}
\newcommand{\cV}{\mathcal V}
\newcommand{\Tr}{\operatorname{Tr}}
\newcommand{\Ad}{\operatorname{Ad}}

\newcommand{\Span}{\operatorname{span}}
\newcommand{\dist}{\operatorname{dist}}
\newcommand{\HS}{\operatorname{HS}}
\newcommand{\op}{\mathrm{op}}
\newcommand{\eps}{\varepsilon}
\newcommand{\E}{\mathrm E}

\begin{document}

\title[Failure of the local lifting property for $C^*(\F_2\times\F_2)$]
{Failure of the local lifting property for the \\ full group C$^*$-algebra of $\F_2\times\F_2$}
\author[A. Ioana]{Adrian Ioana}
\address{Department of Mathematics, University of California San Diego, 9500 Gilman Drive, La Jolla, CA 92093, USA}
\email{aioana@ucsd.edu}

\date{}

\begin{abstract}
We prove that the full group C$^*$-algebra of $\F_2\times\F_2$ does not have the local lifting property. Consequently, neither the lifting property nor the local lifting property is preserved under maximal tensor products. This also gives examples of groups with the Haagerup property whose full group C$^*$-algebras fail the local lifting property. Our proof builds on the asymptotic homomorphisms into matrix algebras constructed in \cite{Ioa21} to prove failure of Hilbert-Schmidt stability for products of free groups. We deduce failure of the local lifting property by showing that these maps cannot even be approximated, in the normalized Hilbert-Schmidt norm, by unital completely positive maps. This result was obtained using ChatGPT Astra.
\end{abstract}
\maketitle

\section{Introduction and statement of the main result}

The local lifting property (LLP), introduced by Kirchberg in \cite{Kir93}, requires that unital completely positive maps from a unital C$^*$-algebra into quotient C$^*$-algebras admit unital completely positive lifts on finite-dimensional operator systems. It is a local version of the lifting property (LP), which requires such maps to admit lifts on the whole algebra. These properties play an important role in the theory of tensor products of C$^*$-algebras and in the study of approximation properties of groups; see \cite{Oz04b,Pis20} for background. In this paper, we study the LLP for full group C$^*$-algebras of countable discrete groups.

The Choi-Effros lifting theorem \cite{CE76} implies that every separable nuclear C$^*$-algebra has the LP. In particular, $C^*(\Gamma)$ has the LP whenever $\Gamma$ is a countable amenable group. Kirchberg proved that full C$^*$-algebras of countable free groups also have the LP \cite[Lemma~3.3]{Kir94}. Recently, Fournier-Facio and Willett enlarged the class of groups known to have this property. Their results apply, in particular, to finitely generated fundamental groups of manifolds of dimension at most three, and hence to surface groups \cite[Theorem~1.2]{FFW26}.

Beyond group algebras, Courtney proved the LLP for certain universal algebras generated by algebraic contractions and the LP for universal Pythagorean algebras \cite[Corollary~3.8 and Proposition~4.1]{Cou21}. Enders and Shulman proved that the LP is preserved under full free products of separable C$^*$-algebras amalgamated over finite-dimensional subalgebras and that Exel's soft tori have the LP \cite[Corollary~4.3 and Theorem~5.5]{ES24}.

Junge and Pisier \cite{JP95} proved that the minimal and maximal tensor norms on $\mathbb B(\ell^2)\odot\mathbb B(\ell^2)$ are distinct. By Kirchberg's tensor product characterization, this means that $\mathbb B(\ell^2)$ does not have the LLP.
For full group C$^*$-algebras, Ozawa showed that there exist countable groups whose full C$^*$-algebras do not have the LP \cite[Corollary~5]{Oz04a}. Thom subsequently constructed an explicit example of a non-residually finite property~(T) group whose full C$^*$-algebra fails the LLP \cite{Tho10}. Further examples were obtained by Buss, Echterhoff and Willett \cite[Corollary~4.8]{BEW18}. In \cite{ISW20}, Spaas, Wiersma and the author used cohomological methods to establish failure of lifting properties for large classes of groups. In particular, we showed that $C^*(\mathbb Z^2\rtimes\operatorname{SL}_2(\mathbb Z))$ and $C^*(\operatorname{SL}_n(\mathbb Z))$, for $n\geq3$, do not have the LLP, and that the full C$^*$-algebra of any infinitely presented property~(T) group fails the LP.

The algebra $C^*(\F_2\times\F_2)$ has played a central role in the theory of tensor products of C$^*$-algebras since Kirchberg's work \cite{Kir93}. Kirchberg showed that Connes' embedding problem is equivalent both to the equality of the minimal and maximal tensor norms on $C^*(\F_2)\odot C^*(\F_2)$ and to residual finite-dimensionality of $C^*(\F_2\times\F_2)$; see also \cite[Theorem~1.1 and Proposition~3.19]{Oz04b}. Thus the negative solution of Connes' embedding problem obtained by Ji, Natarajan, Vidick, Wright and Yuen through their theorem $\mathrm{MIP}^*=\mathrm{RE}$ \cite{JNVWY20} implies that
\[
C^*(\F_2)\otimes_{\min}C^*(\F_2)
\neq C^*(\F_2)\otimes_{\max}C^*(\F_2),
\]
and that $C^*(\F_2\times\F_2)$ is not residually finite-dimensional.

Ozawa also related the LLP for $C^*(\F_2\times\F_2)$ to Tsirelson's problem. His Strong Kirchberg Conjecture is equivalent to the conjunction of Kirchberg's conjecture and the LLP for $C^*(\F_2\times\F_2)$ \cite[Theorem~4]{Oz13}. In its dilation formulation, whenever the commutators between two finite families of unitary matrices are sufficiently small in operator norm, the families can be approximated by compressions of two unitary families that commute with each other. Requiring finite-dimensional dilation spaces gives the Strong Kirchberg Conjecture, while allowing arbitrary Hilbert spaces characterizes the LLP for $C^*(\F_2\times\F_2)$ \cite[Lemma~5]{Oz13}. The negative solution of Connes' embedding problem therefore disproved the Strong Kirchberg Conjecture, but left the LLP question open.

The question of whether $C^*(\F_2\times\F_2)$ has the LLP goes back to \cite[Section~3]{Oz04b} and is discussed in \cite[p.~435]{Pis20}, \cite{Ioa21,ES24} and \cite[Question~1.6]{FFW26}. The main goal of this paper is to give a negative answer.

\begin{main}\label{main}
The full group C$^*$-algebra $C^*(\F_2\times\F_2)$ does not have the LLP.
\end{main}

Theorem~\ref{main} was also obtained independently by Roydor \cite{Roy26}. We discuss the relation between the two proofs below.

Before discussing the proof, we explain several consequences of Theorem~\ref{main}.

Since $C^*(\F_2)$ has the LP and
$C^*(\F_2\times\F_2)\cong C^*(\F_2)\otimes_{\max}C^*(\F_2)$,
Theorem~\ref{main} shows that neither the LP nor the LLP is preserved under maximal tensor products. Thus it also answers the permanence questions recorded in \cite[Section~3]{Oz04b}.

By Kirchberg's tensor product characterization of the LLP \cite{Kir93} (see also \cite[Corollary~3.17(ii)]{Oz04b}), Theorem~\ref{main} gives
\[
C^*(\F_2\times\F_2)\otimes_{\min}\mathbb B(\ell^2)
\neq C^*(\F_2\times\F_2)\otimes_{\max}\mathbb B(\ell^2).
\]
In other words, the minimal and maximal tensor norms on $C^*(\F_2\times\F_2)\odot\mathbb B(\ell^2)$ are distinct.

Enders and Shulman proved that $C^*(\F_2\times\F_2)$ has the LLP if and only if its extension semigroup is a group \cite[Corollary~7.6]{ES24}. Hence Theorem~\ref{main} also implies that $\operatorname{Ext}(C^*(\F_2\times\F_2))$ is not a group.

The group $\F_2\times\F_2$ has the Haagerup property. Theorem~\ref{main} therefore gives examples of groups with the Haagerup property whose full C$^*$-algebras fail the LLP, answering \cite[Question~1.8]{FFW26}. Such groups do not have relative property~(T) with respect to any infinite subgroup. This differs from the obstructions in \cite{Tho10,ISW20}, which use property~(T) or its relative version.

The LLP passes from $C^*(\Gamma)$ to $C^*(\Lambda)$ whenever $\Lambda<\Gamma$; see \cite[Remark~1.2]{ISW20}. We thus obtain the following consequence of Theorem~\ref{main}.

\begin{mcor}\label{products}
Let $\Gamma_1$ and $\Gamma_2$ be countable groups which contain non-abelian free subgroups. Then $C^*(\Gamma_1\times\Gamma_2)$ does not have the LLP. In particular, $C^*(\F_m\times\F_n)$ fails the LLP for all $2\leq m,n\leq\infty$.
\end{mcor}

More generally, Theorem~\ref{main} rules out the LLP for the full C$^*$-algebra of any countable group containing $\F_2\times\F_2$. For example, this applies to every right-angled Artin group whose defining graph contains an induced cycle of length four.
This complements the positive result for chordal defining graphs in \cite[Theorem~1.2]{FFW26}; whether the absence of induced cycles of length four suffices for the LLP remains open \cite[Question~1.7]{FFW26}.

For every $n\geq4$, the arithmetic group $\operatorname{SL}_n(\mathbb Z)$ contains a copy of $\F_2\times\F_2$. Thus Theorem~\ref{main} also recovers the failure of the LLP for $C^*(\operatorname{SL}_n(\mathbb Z))$ for these values of $n$, established in \cite{ISW20}.

\subsection*{Comments on the proof of Theorem~\ref{main}}
For full group C$^*$-algebras, the LLP imposes a useful restriction on asymptotic homomorphisms into matrix algebras. Namely, if $C^*(\Gamma)$ has the LLP, then every sequence of maps $\pi_k:\Gamma\rightarrow\cU(\mathbb M_{d_k}(\C))$ satisfying
\[
\lim_{k\rightarrow\infty}\|\pi_k(gh)-\pi_k(g)\pi_k(h)\|_2=0,
\qquad\text{for every }g,h\in\Gamma,
\]
can be approximated pointwise in $\|\cdot\|_2$ by unital completely positive maps $C^*(\Gamma)\rightarrow\mathbb M_{d_k}(\C)$ \cite[Corollary~1.7]{ISW20}. Here $\|\cdot\|_2$ denotes the normalized Hilbert-Schmidt norm.

In \cite{Ioa21}, we proved that products of two non-abelian free groups are not Hilbert-Schmidt stable, by constructing asymptotic homomorphisms into matrix algebras which cannot be approximated by homomorphisms. The introduction of \cite{Ioa21} raised the question of whether the associated map into a tracial ultraproduct admits a unital completely positive lift. To prove Theorem~\ref{main}, we strengthen the obstruction furnished by that construction: we show that suitable asymptotic homomorphisms cannot be approximated even by unital completely positive maps. Their multiplicative defects in fact tend to zero in the operator norm, while the obstruction to approximation holds in the normalized Hilbert-Schmidt norm.

We start with the finite sets of almost commuting unitaries from \cite[Section~6]{Ioa21}. As in \cite{Ioa21}, quantum expanders \cite{Has07,Pis12} allow us to encode these sets into two triples of unitaries, giving asymptotic homomorphisms of $\F_3\times\F_3$. The key step is Lemma~\ref{state}: assuming the existence of unital completely positive approximations, Stinespring's dilation theorem and the standard form of a von Neumann algebra provide positive Hilbert-Schmidt operators with prescribed commutation and support estimates.

The second free factor, together with the deformation from \cite{Ioa21}, forces these operators to be concentrated on Walsh functions of low degree. The first free factor forces them to almost commute, on average, with all coordinate functions. The Fourier estimate in Lemma~\ref{fourier} shows that these two conditions are incompatible when the degree is small compared with the number of coordinates. Finally, the inclusion $\F_3\times\F_3<\F_2\times\F_2$ yields Theorem~\ref{main}.

The paper is organized as follows. In Section~2, we recall the required facts about the LLP, quantum expanders and standard forms, and prove the two lemmas used in the obstruction argument. Section~3 contains the finite matrix construction, including proofs of the deformation identities from \cite{Ioa21} which we need. We prove Theorem~\ref{main} in Section~4.

\subsection*{Independent work}
While preparing this manuscript for submission to arXiv, the author learned of Roydor's preprint \cite{Roy26}, which independently proves Theorem~\ref{main}.

Both proofs build on the almost commuting matrix construction and quantum expander arguments of \cite{Ioa21}, and use standard forms to obtain positive vectors from almost intertwining relations; compare Lemma~\ref{state} with \cite[Lemma~3.1]{Roy26}. Roydor realizes the deformation using rotated Walsh bases in a single matrix algebra and proves a uniform estimate for almost intertwiners \cite[Theorem~D]{Roy26}. His concluding estimate combines control of the Walsh degree with a matrix spectral gap and a unitary exchanging complementary subspaces \cite[Proposition~4.2]{Roy26}. Here we retain the doubled tensor product deformation from \cite{Ioa21} and derive the contradiction from the Fourier estimate in Lemma~\ref{fourier}. Roydor also proves that $C^*(\F_2\times\F_2)$ does not locally embed into $C^*(\F_\infty)$ \cite[Theorem~C]{Roy26}.

\subsection*{Acknowledgements}
The author was supported in part by the NSF grant DMS-2451697 and the Presidential Chair in Mathematics. He also thanks OpenAI for providing access to ChatGPT through the ChatGPT for Academic Researchers program.

\subsection*{AI declaration}
The initial proposed proof of Theorem~\ref{main} was generated by ChatGPT (Astra) on November 15, 2026 (Pacific Standard Time), in response to the author's request to prove that $C^*(\F_2\times\F_2)$ does not have the local lifting property. In that request, the author suggested using the failure of Hilbert-Schmidt stability established in \cite{Ioa21} as a possible approach. The author had previously tried this approach without success while working on \cite{Ioa21} and had recorded the resulting question in its introduction. ChatGPT supplied the crucial new Lemma~\ref{state} and the Fourier estimate in Lemma~\ref{fourier}. The latter addresses, in the setting of Hilbert-Schmidt operators, the obstruction treated for subalgebras by the dimension estimate in \cite[Lemma~6.6]{Ioa21}.

The resulting \LaTeX{} draft served as the starting point for this manuscript. In subsequent interactions, AI assistance was used to expand and refine the arguments, organize the finite matrix construction and the preliminaries, locate related literature, suggest consequences of the main result, and draft and revise the exposition. The author subsequently revised the exposition and the preliminaries and selected the consequences retained in the manuscript. AI assistance was also used to proofread the manuscript, check the related literature, clarify the Fourier estimate, and improve the notation and presentation. The author takes responsibility for the mathematical content, the citations, and the final presentation.

\section{Preliminaries}

\subsection{Notation}
For a tracial von Neumann algebra $(M,\tau)$, we denote by $\|x\|_2=\tau(x^*x)^{1/2}$ its $2$-norm and by $L^2(M)$ the associated Hilbert space. We use the inner product $\langle x,y\rangle=\tau(y^*x)$. If $Q\subset M$ is a unital von Neumann subalgebra, then $\E_Q:M\rightarrow Q$ denotes the trace-preserving conditional expectation. We use the same notation for its extension to the orthogonal projection from $L^2(M)$ onto $L^2(Q)$. We write $\cU(M)$ for the unitary group of $M$ and $[x,y]=xy-yx$.

Every full matrix algebra $\mathbb M_d(\C)$ is endowed with its normalized trace $\tau_d=d^{-1}\Tr$. For a Hilbert space $\cH$, we denote by $\HS(\cH)$ the Hilbert space of Hilbert-Schmidt operators on $\cH$, with norm $\|T\|_{\HS}=\Tr(T^*T)^{1/2}$. This norm and the trace norm $\|T\|_1=\Tr(|T|)$ use the ordinary, unnormalized trace. Keeping these two normalizations separate will be useful below.

For a matrix algebra $N$, put $K_N=L^2(N)$ and denote by $L_a$ and $R_a$ the operators $L_a(x)=ax$ and $R_a(x)=xa$ on $K_N$. Thus $a^{\op}\mapsto R_a$ is a representation of $N^{\op}$, since $R_aR_b=R_{ba}$. For $b\in\cU(N)$, we write $\Ad(b)=L_bR_{b^*}$ for the unitary operator $x\mapsto bxb^*$ on $K_N$.

\subsection{The local lifting property}
An operator system in a unital C$^*$-algebra is a self-adjoint linear subspace containing the unit. We abbreviate ``unital completely positive'' by u.c.p.

\begin{definition}
A unital C$^*$-algebra $A$ has the \emph{local lifting property} if, whenever $B$ is a unital C$^*$-algebra, $J\subset B$ is a closed two-sided ideal and $\varphi:A\rightarrow B/J$ is a u.c.p. map, the following holds: for every finite-dimensional operator system $E\subset A$, there is a u.c.p. map $\psi:E\rightarrow B$ such that $q\circ\psi=\varphi|_E$, where $q:B\rightarrow B/J$ is the quotient map. If every such map $\varphi$ admits a u.c.p. lift $\psi:A\rightarrow B$, then $A$ has the \emph{lifting property} (LP).
\end{definition}

For a countable group $\Gamma$, we denote by $\{u_g\mid g\in\Gamma\}$ the canonical unitaries in $C^*(\Gamma)$. A sequence of maps $\pi_k:\Gamma\rightarrow\cU(N_k)$, where $(N_k,\tau_k)$ are tracial von Neumann algebras, is an \emph{asymptotic homomorphism} if $\lim_{k\rightarrow\infty}\|\pi_k(gh)-\pi_k(g)\pi_k(h)\|_{2,\tau_k}=0$, for every $g,h\in\Gamma$.

We recall the following result, which in particular applies when the $N_k$ are matrix algebras. Here, $R$ denotes the hyperfinite II$_1$ factor and $\omega$ is a free ultrafilter on $\N$.

\begin{corollary}[{\cite[Corollary~1.7]{ISW20}}]\label{ISW}
Let $\Gamma$ be a countable group, let $(N_k,\tau_k)$ be tracial von Neumann algebras, and let $\pi_k:\Gamma\rightarrow\cU(N_k)$ be an asymptotic homomorphism. Assume either that $C^*(\Gamma)$ has the LP, or that $C^*(\Gamma)$ has the LLP and every $(N_k,\tau_k)$ admits a trace-preserving embedding into $R^\omega$. Then there are u.c.p. maps $\varphi_k:C^*(\Gamma)\rightarrow N_k$ such that
\[
\lim_{k\rightarrow\infty}\|\varphi_k(u_g)-\pi_k(g)\|_{2,\tau_k}=0,
\qquad\text{for every }g\in\Gamma.
\]
\end{corollary}

\subsection{Quantum expanders}
We will use the following form of the quantum expander theorem of Hastings \cite{Has07}; see also \cite[Lemma~1.8]{Pis12} and \cite[Lemma~4.3]{Ioa21}.

\begin{proposition}\label{gap}
There is a constant $\kappa>0$ such that for every $d\in\N$, there are unitaries $c_1,c_2\in\mathbb M_d(\C)$ satisfying
\begin{equation}\label{gap-matrix}
\|x-\tau_d(x)1\|_2^2\leq\kappa^2\sum_{i=1}^2\|[c_i,x]\|_2^2,
\qquad x\in\mathbb M_d(\C).
\end{equation}
The following two consequences hold with the same constant:
\begin{enumerate}[label=(\arabic*)]
\item If $\cL$ is a Hilbert space and $e_0$ is the projection onto $\C1\subset L^2(\mathbb M_d(\C))$, then
\begin{equation}\label{gap-vectors}
\|(1\otimes(1-e_0))\eta\|^2
\leq\kappa^2\sum_{i=1}^2
\|(1\otimes\Ad(c_i))\eta-\eta\|^2
\end{equation}
for every $\eta\in\cL\otimes L^2(\mathbb M_d(\C))$.
\item Let $\sigma:\mathbb M_d(\C)\rightarrow\mathbb B(\cH)$ be a unital representation and let $\mathcal E$ be the orthogonal projection of $\HS(\cH)$ onto the operators commuting with $\sigma(\mathbb M_d(\C))$. Then
\begin{equation}\label{gap-operators}
\|T-\mathcal E(T)\|_{\HS}^2
\leq\kappa^2\sum_{i=1}^2\|[T,\sigma(c_i)]\|_{\HS}^2,
\qquad T\in\HS(\cH).
\end{equation}
The analogous statement holds for representations of $\mathbb M_d(\C)^{\op}$ with the pair $c_1^{\op},c_2^{\op}$.
\end{enumerate}
\end{proposition}

\begin{proof}
The main assertion \eqref{gap-matrix} follows immediately from \cite[Lemma~4.3]{Ioa21}.
%For a triple $v=(v_1,v_2,v_3)\in\cU(\mathbb M_d(\C))^3$, let $T_v$ be the operator on $L^2(\mathbb M_d(\C))$ given by $T_v(x)=\sum_{i=1}^3v_ixv_i^*$, and denote its restriction to $L^2(\mathbb M_d(\C))\ominus\C1$ by $T_v^0$. Hastings' theorem, in the formulation of \cite[Lemma~1.8]{Pis12}, gives a fixed $\alpha<3$ such that, for every sufficiently large $d$, there is a triple $v$ with $\|T_v^0\|\leq\alpha$. Indeed, one can take any $\alpha$ strictly between $2\sqrt2$ and $3$.

%Put $c_1=v_3^*v_1$ and $c_2=v_3^*v_2$. Composing $T_v$ with $\Ad(v_3^*)$ shows that the restriction of $\id+\Ad(c_1)+\Ad(c_2)$ to the trace-zero subspace has norm at most $\alpha$. Consequently, if $\tau_d(x)=0$, then \[ (3-\alpha)\|x\|_2 \leq\|x-c_1xc_1^*\|_2+\|x-c_2xc_2^*\|_2. \]

%Squaring and applying Cauchy--Schwarz gives \eqref{gap-matrix} for trace-zero $x$, with $\kappa=\sqrt2/(3-\alpha)$. Applying this to $x-\tau_d(x)1$ proves \eqref{gap-matrix} for arbitrary $x$ in these dimensions. The case $d=1$ is immediate. In each of the finitely many remaining dimensions $d\geq2$, choose a pair generating the full matrix algebra. The only matrices commuting with both unitaries are scalar, so compactness of the unit sphere in the trace-zero subspace gives a finite constant in \eqref{gap-matrix}. Increasing $\kappa$ proves the assertion for every $d$.

To obtain \eqref{gap-vectors}, expand $\eta$ in an orthonormal basis of $\cL$, apply \eqref{gap-matrix} to each coefficient, and add the resulting inequalities. To prove \eqref{gap-operators}, identify $\cH=\C^d\otimes\cL$ so that $\sigma(x)=x\otimes1$. Under
\[
\HS(\cH)=\HS(\C^d)\otimes\HS(\cL),
\]
conjugation by $\sigma(c_i)$ is $\Ad(c_i)\otimes1$ and $\mathcal E$ projects the first factor onto $\C I_d$. Multiplying \eqref{gap-matrix} by $d$ converts it to the ordinary Hilbert-Schmidt normalization, so the same coefficient argument applies. Equivalently,
\[
\mathcal E(T)=\int_{\cU(\mathbb M_d(\C))}\sigma(u)T\sigma(u)^*\,du.
\]
Finally, %the transpose identification of $\mathbb M_d(\C)^{\op}$ with $\mathbb M_d(\C)$ preserves the trace and the $2$-norm, and sends commutators to their negatives after transposition. Thus 
it is easy to see that the same bound applies to the opposite algebra.
\end{proof}

\subsection{Standard forms and positive vectors}
We recall the standard form facts used below; see \cite{Haa75} and \cite[Section~2.2]{AHW13}. Every von Neumann algebra $P$ has a standard form on $L^2(P)$ with a distinguished positive cone. Each normal positive functional $\omega$ on $P$ is represented by a unique vector $\xi_\omega$ in that cone, and
\begin{equation}\label{cone}
\|\xi_\omega-\xi_\nu\|^2\leq\|\omega-\nu\|,
\qquad \omega,\nu\in P_*^+.
\end{equation}
Automorphisms have canonical unitary implementations preserving the cone. For a full matrix algebra $N$, the standard Hilbert space of $P\overline\otimes N^{\op}$ can be identified with $L^2(P)\otimes L^2(N)$, with $p\otimes x^{\op}$ acting as $p\otimes R_x$. The algebra $P$ need not be finite or admit a trace.

We will also use the Powers-St\o rmer inequality \cite{PS70} in the form
\begin{equation}\label{PS}
\|S^{1/2}-T^{1/2}\|_{\HS}^2\leq\|S-T\|_1,
\qquad S,T\geq0\text{ trace class}.
\end{equation}
For vectors $\eta,\zeta$ in a Hilbert space, let $r_{\eta,\zeta}$ denote the rank-one operator
\[
r_{\eta,\zeta}(v)=\langle v,\zeta\rangle\eta.
\]
Thus $\|r_{\eta,\zeta}\|_1=\|\eta\|\,\|\zeta\|$. For unit vectors $\eta,\zeta$, their vector states and density operators satisfy
\begin{equation}\label{rank-one}
\|\omega_\eta-\omega_\zeta\|
\leq\|r_{\eta,\eta}-r_{\zeta,\zeta}\|_1
\leq2\|\eta-\zeta\|.
\end{equation}
The last inequality follows from
$r_{\eta,\eta}-r_{\zeta,\zeta}=r_{\eta-\zeta,\eta}+r_{\zeta,\eta-\zeta}$.

\subsection{A consequence of completely positive approximation}
The next lemma allows us to treat u.c.p. maps instead of homomorphisms.

\begin{lemma}\label{state}
Let $u_1,\ldots,u_r$ and $v_1,\ldots,v_s$ be the canonical generators of $C^*(\F_r\times\F_s)$, and let $N$ be a full matrix algebra. Suppose that $a_i,b_j\in\cU(N)$ and a u.c.p. map $\varphi:C^*(\F_r\times\F_s)\rightarrow N$ satisfy
\[
\|\varphi(u_i)-a_i\|_2\leq\eps,
\qquad \|\varphi(v_j)-b_j\|_2\leq\eps,
\quad 1\leq i\leq r,\ 1\leq j\leq s.
\]
Then there are a Hilbert space $\cK$, a unit vector $\xi\in\cK\otimes K_N$ and a positive operator $h$ on $K_N$ with $\|h\|_{\HS}=1$ such that, with $c=2^{3/4}$,
\begin{align}
\|(1\otimes\Ad(b_j))\xi-\xi\|&\leq c\eps^{1/4},\label{state-b}\\
\|[h,R_{a_i}]\|_{\HS}&\leq c\eps^{1/4},\label{state-a}\\
\|(1-e)h\|_{\HS}&=\|(1\otimes(1-e))\xi\|\label{state-support}
\end{align}
for every $i,j$ and every orthogonal projection $e$ on $K_N$.
\end{lemma}

\begin{proof}
Write $N=\mathbb M_d(\C)$. By Stinespring's dilation theorem \cite{Sti55}, there are a unital representation $\pi:C^*(\F_r\times\F_s)\rightarrow\mathbb B(H)$ and an isometry $V:\C^d\rightarrow H$ such that $\varphi(x)=V^*\pi(x)V$. Put
\[
p_i=\pi(u_i),\qquad q_j=\pi(v_j),\qquad
P=\{p_1,\ldots,p_r\}'',\qquad \delta=(2\eps)^{1/2}.
\]
The defining relations of the product group imply that $q_j\in P'$, for every $j$.

\smallskip
\noindent\emph{Step 1. Estimates for the Stinespring isometry.}
For a linear operator $W:\C^d\rightarrow H$, we set $\|W\|_{2,d}^2=d^{-1}\Tr(W^*W)$. Since $p_i$ and $a_i$ are unitaries and $V^*V=1$, we have
\begin{align*}
\|p_iV-Va_i\|_{2,d}^2
&=2-2\operatorname{Re}\tau_d(a_i^*\varphi(u_i))\\
&=2\operatorname{Re}\tau_d\big(a_i^*(a_i-\varphi(u_i))\big)
\leq2\|\varphi(u_i)-a_i\|_2\leq2\eps.
\end{align*}
The same argument gives
\begin{equation}\label{isometry}
\|p_iV-Va_i\|_{2,d}\leq\delta,
\qquad \|q_jV-Vb_j\|_{2,d}\leq\delta.
\end{equation}
Multiplying the second difference on the left by $q_j^*$ and on the right by $b_j^*$ also gives
\begin{equation}\label{isometry-adjoint}
\|q_j^*V-Vb_j^*\|_{2,d}\leq\delta.
\end{equation}

\smallskip
\noindent\emph{Step 2. A normal state and its positive vector.}
Represent $P\overline\otimes N^{\op}$ normally on $\HS(\C^d,H)$ by
\[
(p\otimes x^{\op})W=pWx.
\]
Since $N$ is finite-dimensional, the commuting left and right actions indeed define a normal representation of this spatial tensor product. The vector $V/\sqrt d$ has norm one in the ordinary Hilbert-Schmidt norm. Its vector state is
\begin{equation}\label{omega}
\omega(p\otimes x^{\op})=\tau_d(V^*pVx).
\end{equation}
This is a normal state, with no finiteness assumption on $P$.

Let $\xi$ be the unique positive-cone vector for $\omega$ in the standard Hilbert space $L^2(P)\otimes K_N$, and put $\cK=L^2(P)$. Then $\|\xi\|=1$. Both $V/\sqrt d$ and $\xi$ implement the same state. Applying that state to the positive element
\[
(p_i\otimes1-1\otimes a_i^{\op})^*
(p_i\otimes1-1\otimes a_i^{\op})
\]
and using \eqref{isometry}, we obtain
\begin{equation}\label{state-intertwining}
\|(p_i\otimes1)\xi-(1\otimes R_{a_i})\xi\|\leq\delta.
\end{equation}
Here only equality of the two vector states is used.

\smallskip
\noindent\emph{Step 3. Approximate invariance under the second factor.}
Fix $j$ and write $q=q_j$, $b=b_j$. Since $q$ commutes with $P$, the vector $q^*V/\sqrt d$ defines the same state $\omega$ as $V/\sqrt d$. Indeed,
\[
\tau_d\big((q^*V)^*p(q^*V)x\big)
=\tau_d(V^*q p q^*Vx)=\tau_d(V^*pVx).
\]
By \eqref{isometry-adjoint}, this vector is within $\delta$ of $Vb^*/\sqrt d$. The state implemented by the latter vector is
\[
\omega_b=\omega\circ\Ad(1\otimes b^{\op}).
\]
To check the order of the factors, observe that
$b^{\op}x^{\op}(b^*)^{\op}=(b^*xb)^{\op}$, and use cyclicity of $\tau_d$. Thus \eqref{rank-one} gives
\begin{equation}\label{state-distance}
\|\omega_b-\omega\|\leq2\delta.
\end{equation}

On $K_N$, let $J_N(x)=x^*$. The standard implementation of the inner automorphism $\Ad(b^{\op})$ is
\[
R_bJ_NR_bJ_N=L_{b^*}R_b=\Ad(b^*).
\]
Consequently the positive-cone vector of $\omega_b=\omega\circ\Ad(1\otimes b^{\op})$ is $(1\otimes\Ad(b))\xi$. Combining \eqref{cone} and \eqref{state-distance}, we get
\[
\|(1\otimes\Ad(b))\xi-\xi\|^2\leq2\delta.
\]
Since $(2\delta)^{1/2}=2^{3/4}\eps^{1/4}$, this proves \eqref{state-b}.

\smallskip
\noindent\emph{Step 4. The reduced density operator.}
Let $\overline{\cK}$ be the conjugate Hilbert space of $\cK$, with $\overline p(\overline\eta)=\overline{p\eta}$. We use the conjugate-linear isometric identification
\[
\cK\otimes K_N\longrightarrow\HS(K_N,\overline{\cK}),
\qquad
\eta\otimes x\longmapsto\big(y\mapsto\langle y,x\rangle\overline\eta\big).
\]
We also write $\xi$ for the corresponding operator $K_N\rightarrow\overline{\cK}$, with $\|\xi\|_{\HS}=1$. Then $(p\otimes z)\xi$ corresponds to $\overline p\,\xi z^*$, for $p\in P$ and $z\in\mathbb B(K_N)$. Define
\begin{equation}\label{density-operator}
T=\xi^*\xi\quad\text{on }K_N,
\qquad h=T^{1/2}.
\end{equation}
Then $T\geq0$, $\Tr(T)=\|\xi\|_{\HS}^2=1$, and $\|h\|_{\HS}=1$.

By \eqref{state-intertwining}, we have
\[
\|\overline{p_i}\xi-\xi R_{a_i}^*\|_{\HS}\leq\delta.
\]
Put $A_i=\overline{p_i}\xi$ and $B_i=\xi R_{a_i}^*$. Since $\overline{p_i}$ and $R_{a_i}$ are unitaries, $\|A_i\|_{\HS}=\|B_i\|_{\HS}=1$ and
\[
A_i^*A_i=T,
\qquad B_i^*B_i=R_{a_i}TR_{a_i}^*.
\]
Using the inequality $\|X^*Y\|_1\leq\|X\|_{\HS}\|Y\|_{\HS}$ for Hilbert-Schmidt operators, we obtain
\begin{equation}\label{density-invariance}
\begin{aligned}
\|T-R_{a_i}TR_{a_i}^*\|_1
&=\|A_i^*A_i-B_i^*B_i\|_1\\
&\leq\|A_i^*(A_i-B_i)\|_1
+\|(A_i^*-B_i^*)B_i\|_1\\
&\leq(\|A_i\|_{\HS}+\|B_i\|_{\HS})\|A_i-B_i\|_{\HS}
\leq2\delta.
\end{aligned}
\end{equation}
Applying \eqref{PS} and using $(R_{a_i}TR_{a_i}^*)^{1/2}=R_{a_i}hR_{a_i}^*$ gives
\[
\|[h,R_{a_i}]\|_{\HS}^2
=\|h-R_{a_i}hR_{a_i}^*\|_{\HS}^2
\leq2\delta.
\]
Since $(2\delta)^{1/2}=2^{3/4}\eps^{1/4}$, this proves \eqref{state-a}.

Finally, for an orthogonal projection $e$ on $K_N$, the vector $(1\otimes(1-e))\xi$ corresponds to the operator $\xi(1-e)$. By cyclicity of the trace and \eqref{density-operator},
\begin{align*}
\|(1-e)h\|_{\HS}^2
&=\Tr(h(1-e)h)=\Tr((1-e)T)\\
&=\|\xi(1-e)\|_{\HS}^2
=\|(1\otimes(1-e))\xi\|^2.
\end{align*}
This proves \eqref{state-support} and finishes the proof.
\end{proof}

\subsection{A Fourier estimate}
Let $n\geq1$ be an integer and equip $\{-1,1\}^n$ with the uniform probability measure. For $S\subset\{1,\ldots,n\}$, put $w_S(x)=\prod_{i\in S}x_i$. These functions form an orthonormal basis of $L^2(\{-1,1\}^n)$. Let $M_i$ be multiplication by the $i$-th coordinate function, so that $M_iw_S=w_{S\triangle\{i\}}$. For an integer $0\leq l\leq n$, denote by $f_l$ the orthogonal projection onto $\Span\{w_S\mid |S|\leq l\}$.

\begin{lemma}\label{fourier}
Let $\cL$ be a Hilbert space and $y\in\HS(\cL\otimes L^2(\{-1,1\}^n))$. If $y=(1\otimes f_l)y(1\otimes f_l)$, then
\begin{equation}\label{fourier-bound}
\frac1n\sum_{i=1}^n\|[y,1\otimes M_i]\|_{\HS}^2
\geq2\left(1-2\sqrt{\frac ln}\right)\|y\|_{\HS}^2.
\end{equation}
\end{lemma}

\begin{proof}
Put $\cH=L^2(\{-1,1\}^n)$. It suffices to prove the inequality when $\cL=\C$.

% Indeed, relative to an orthonormal basis of $\cL$, write $y$ as a block matrix $(y_{\alpha,\beta})$, with $y_{\alpha,\beta}\in\mathbb B(\cH)$. Then $y_{\alpha,\beta}=f_ly_{\alpha,\beta}f_l$ and \[ \|y\|_{\HS}^2=\sum_{\alpha,\beta}\|y_{\alpha,\beta}\|_{\HS}^2, \qquad \|[y,1\otimes M_i]\|_{\HS}^2 =\sum_{\alpha,\beta}\|[y_{\alpha,\beta},M_i]\|_{\HS}^2. \] Thus applying the scalar case to each block and summing gives the general case.

Assuming that $\cL=\C$, write $y=(y_{S,T})$ in the Walsh basis. Put
\[
\mathscr I_l=\{(S,T)\mid |S|,|T|\leq l\},
\qquad z=(y_{S,T})_{(S,T)\in\mathscr I_l}\in\ell^2(\mathscr I_l).
\]
The support assumption gives $y_{S,T}=0$ outside $\mathscr I_l$, so $\|z\|^2=\|y\|_{\HS}^2$. Conjugation by $M_i$ changes the coefficient index $(S,T)$ to $(S\triangle\{i\},T\triangle\{i\})$. Let $A$ be the adjacency operator of the graph on $\mathscr I_l$ joining these two indices whenever both belong to $\mathscr I_l$. Then
\[
\sum_{i=1}^n\langle M_iyM_i,y\rangle_{\HS}=\langle Az,z\rangle.
\]

Orient each edge from its endpoint with smaller $|S|$ to its endpoint with larger $|S|$, and define $B$ on $\ell^2(\mathscr I_l)$ by
\[
(Bz)_w=\sum_{v\rightarrow w}z_v.
\]
Each vertex has at most $n$ outgoing edges. A vertex $(S,T)$ has at most $|S|\leq l$ incoming edges. 
%since each predecessor is obtained by removing an element of $S$ and toggling the same element in $T$. 
Hence Cauchy--Schwarz gives
\[
\|Bz\|^2
=\sum_w\left|\sum_{v\rightarrow w}z_v\right|^2
\leq l\sum_w\sum_{v\rightarrow w}|z_v|^2
\leq nl\|z\|^2.
\]
Since $A=B+B^*$, we get $\|A\|\leq2\sqrt{nl}$ and therefore
\begin{equation}\label{fourier-adjacency}
\left|\sum_{i=1}^n\langle M_iyM_i,y\rangle_{\HS}\right|
=|\langle Az,z\rangle|
\leq2\sqrt{nl}\|y\|_{\HS}^2.
\end{equation}
Finally, since each $M_i$ is a self-adjoint unitary, \eqref{fourier-adjacency} implies
\begin{align*}
\sum_{i=1}^n\|[y,M_i]\|_{\HS}^2
&=2n\|y\|_{\HS}^2
-2\operatorname{Re}\sum_{i=1}^n\langle M_iyM_i,y\rangle_{\HS}\\
&\geq2(n-2\sqrt{nl})\|y\|_{\HS}^2.
\end{align*}
Dividing by $n$ proves \eqref{fourier-bound}.
\end{proof}

\section{The finite matrix construction}\label{construction-section}

\subsection{Notation and deformation}
We use the notation of \cite[Notation 6.1]{Ioa21}. We include the needed deformation identities and their proofs so that the construction is self-contained.

\begin{notation}\label{construction}
Let $n\in\N$ and $t\in\R$.
\begin{enumerate}[label=(\arabic*)]
\item Put
\[
M_n=\bigotimes_{k=1}^n\mathbb M_2(\C),
\qquad A_n=\bigotimes_{k=1}^n\C^2\subset M_n,
\]
where $\C^2\subset\mathbb M_2(\C)$ is the diagonal subalgebra. All tensor products carry the product normalized trace.
\item Let $\sigma=\begin{pmatrix}1&0\\0&-1\end{pmatrix}$ and denote by $X_{n,i}\in\cU(A_n)$ the tensor with $\sigma$ in position $i$ and $1$ in every other position. For $S\subset\{1,\ldots,n\}$, put
\[
w_S=\prod_{i\in S}X_{n,i}\in A_n.
\]
The identification $A_n=L^\infty(\{-1,1\}^n)$ identifies these elements with the Walsh functions above.
\item Choose a finite subgroup $G_n\subset\cU(A_n\otimes M_n)$ whose linear span is $A_n\otimes M_n$. %For example, use the group generated by $X_{n,i}\otimes1$, $1\otimes X_{n,i}$ and the coordinate copies, in the second factor, of $\begin{pmatrix}0&1\\1&0\end{pmatrix}$. These are tensor products of sign matrices and elements of the finite Pauli group, and their linear span is the indicated algebra. 
In particular,
\begin{equation}\label{group-commutant}
G_n'\cap(M_n\otimes M_n)=A_n\otimes1.
\end{equation}
\item Let $P$ be the orthogonal projection of $\C^2\otimes\C^2$ onto the line spanned by $e_1\otimes e_2-e_2\otimes e_1$, and put
\[
U_t=P+e^{it}(1-P),\qquad \rho_t=\frac{1+\cos t}{2}.
\]
Identify $M_n\otimes M_n$ with $\bigotimes_{k=1}^n(\mathbb M_2(\C)\otimes\mathbb M_2(\C))$ by pairing the corresponding tensor coordinates, and define
\[
\theta_{t,n}=\bigotimes_{k=1}^n\Ad(U_t).
\]
\item Set $\cU_n=\{X_{n,i}\otimes1\mid1\leq i\leq n\}$ and $\cV_{t,n}=G_n\cup\theta_{t,n}(G_n)$. When enumerating these unitaries below, we retain both copies of an element which belongs to the intersection.
\end{enumerate}
\end{notation}

The family $(\theta_{t,n})_{t\in\R}$ is a finite-dimensional analogue of Popa's malleable deformation for noncommutative Bernoulli actions \cite{Po03}, as discussed in \cite[Introduction]{Ioa21}. It joins the identity at $t=0$ to the flip of the two tensor factors at $t=\pi$.

As in \cite[Section~6]{Ioa21}, let $e_l$ denote the orthogonal projection on $L^2(M_n)$ onto the span of tensors with at most $l$ non-scalar coordinates. %More precisely, each tensor coordinate is either $1$ or an element of $\mathbb M_2(\C)\ominus\C1$, and at most $l$ coordinates are of the latter type.
 Its restriction to $L^2(A_n)$ is the Walsh projection $f_l$ from Lemma~\ref{fourier}.

\begin{lemma}[{\cite[Lemma~6.2]{Ioa21}}]\label{deformation}
The following hold:
\begin{enumerate}[label=(\arabic*)]
\item $P(x\otimes1)P=\tau_2(x)P$, for every $x\in\mathbb M_2(\C)$.
\item For $x,y\in\mathbb M_2(\C)$,
\begin{equation}\label{deformation-pairing}
\tau_4\big((x\otimes1)U_t(y\otimes1)U_t^*\big)
=\rho_t\tau_2(xy)+(1-\rho_t)\tau_2(x)\tau_2(y).
\end{equation}
\item For $x\in\mathbb M_2(\C)$ with $\tau_2(x)=0$,
\[
\E_{U_t(\mathbb M_2(\C)\otimes1)U_t^*}(x\otimes1)
=\rho_t U_t(x\otimes1)U_t^*.
\]
\item For $x\in M_n$ and $0\leq l\leq n$,
\begin{equation}\label{deformation-length}
\|\E_{\theta_{t,n}(M_n\otimes1)}(x\otimes1)\|_2^2
\leq(1-\rho_t^{2(l+1)})\|e_l(x)\|_2^2
+\rho_t^{2(l+1)}\|x\|_2^2.
\end{equation}
\item For every $S\subset\{1,\ldots,n\}$,
\begin{equation}\label{deformation-walsh}
\E_{\theta_{t,n}(A_n\otimes1)}(w_S\otimes1)
=\rho_t^{|S|}\theta_{t,n}(w_S\otimes1).
\end{equation}
\end{enumerate}
\end{lemma}

\begin{proof}
Let $\zeta=2^{-1/2}(e_1\otimes e_2-e_2\otimes e_1)$. Since $P$ is the projection onto $\C\zeta$ and $\langle(x\otimes1)\zeta,\zeta\rangle=\tau_2(x)$, we obtain (1).

To prove (2), write $X=x\otimes1$, $Y=y\otimes1$ and $U_t=e^{it}1+(1-e^{it})P$. Expanding and taking traces gives
\begin{align*}
\tau_4(XU_tYU_t^*)
={}&\tau_4(XY)+(e^{it}-1)\tau_4(XYP)\\
&+(e^{-it}-1)\tau_4(XPY)
+|1-e^{it}|^2\tau_4(XPYP).
\end{align*}
Part (1), cyclicity of the trace, and $\tau_4(P)=1/4$ imply
\[
\tau_4(XYP)=\tau_4(XPY)=\tfrac14\tau_2(xy),
\quad
\tau_4(XPYP)=\tfrac14\tau_2(x)\tau_2(y).
\]
Substituting these equalities proves \eqref{deformation-pairing}. Part (3) follows by pairing both sides with $U_t(y\otimes1)U_t^*$ and using (2).

Set $e_{-1}=0$ and let $g_q=e_q-e_{q-1}$ be the orthogonal projection onto tensors of length exactly $q$, for $0\leq q\leq n$. Tensoring (3) with the identity on the scalar coordinates gives
\begin{equation}\label{deformation-exact}
\E_{\theta_{t,n}(M_n\otimes1)}(x\otimes1)
=\sum_{q=0}^n\rho_t^q\theta_{t,n}(g_q(x)\otimes1).
\end{equation}
The terms in this sum are orthogonal. Therefore its squared norm is $\sum_{q=0}^n\rho_t^{2q}\|g_q(x)\|_2^2$. For $q>l$, use $\rho_t^{2q}\leq\rho_t^{2(l+1)}$, and for $q\leq l$, use $\rho_t^{2q}\leq1$. This proves (4), a slightly sharper version of the bound in \cite[Lemma~6.2(4)]{Ioa21}.

For $x=w_S$, formula \eqref{deformation-exact} reduces to $\rho_t^{|S|}\theta_{t,n}(w_S\otimes1)$, which belongs to $\theta_{t,n}(A_n\otimes1)$. Applying the conditional expectation onto this smaller algebra proves (5).
\end{proof}

We record explicitly the Hilbert space consequence of (5). If $\cL$ is any Hilbert space and
\[
\eta=\sum_{S\subset\{1,\ldots,n\}}\eta_S\otimes(w_S\otimes1)
\in\cL\otimes L^2(A_n\otimes1),
\]
then
\begin{equation}\label{deformation-distance}
\dist\big(\eta,\cL\otimes L^2(\theta_{t,n}(A_n\otimes1))\big)^2
=\sum_S(1-\rho_t^{2|S|})\|\eta_S\|^2.
\end{equation}
Indeed, conditional expectation is orthogonal projection and the vectors $\theta_{t,n}(w_S\otimes1)$ form an orthonormal basis of the second algebra.

\begin{lemma}[{\cite[Lemma~6.3]{Ioa21}}]\label{almost-commuting}
For $U\in\cU_n$ and $V\in\cV_{t,n}$, we have
\[
\|[U,V]\|_2\leq\|[U,V]\|\leq4|t|.
\]
\end{lemma}

\begin{proof}
If $V\in G_n$, then $[U,V]=0$, since $U\in A_n\otimes1$ and $G_n\subset A_n\otimes M_n$. Otherwise write $V=\theta_{t,n}(g)$, with $g\in G_n$. Then $[\theta_{t,n}(U),V]=0$, so
\begin{align*}
\|[U,V]\|
&\leq2\|U-\theta_{t,n}(U)\|\\
&=2\|\sigma\otimes1-U_t(\sigma\otimes1)U_t^*\|
\leq4\|U_t-1\|
=4|e^{it}-1|\leq4|t|.
\end{align*}
The inequality between the normalized $2$-norm and the operator norm proves the remaining assertion.
\end{proof}

\begin{lemma}\label{averaging}
Let $B$ be a finite-dimensional tracial von Neumann algebra and $G\subset\cU(B)$ a finite subgroup. Put $Q=G'\cap B$. Then
\[
\frac1{|G|}\sum_{g\in G}\Ad(g)=\E_Q
\quad\text{on }L^2(B).
\]
Consequently, for every Hilbert space $\cL$ and $\eta\in\cL\otimes L^2(B)$,
\begin{equation}\label{averaging-distance}
\frac1{|G|}\sum_{g\in G}\|(1\otimes\Ad(g))\eta-\eta\|^2
=2\|(1\otimes(1-\E_Q))\eta\|^2.
\end{equation}
\end{lemma}

\begin{proof}
The average is a self-adjoint idempotent, since $G$ is a group. Its range consists precisely of the vectors fixed by all $\Ad(g)$, namely $L^2(Q)$. This proves the first equality. Expand the squared norms in the second equality and use the first one.
\end{proof}

\subsection{Encoding the finite sets by two triples}
Fix $n\in\N$ and $t\in\R$. Enumerate $G_n=\{g_1,\ldots,g_q\}$, put $m=2q$, and define
\[
V_j=g_j\ (1\leq j\leq q),\qquad
V_{q+j}=\theta_{t,n}(g_j)\ (1\leq j\leq q).
\]
This list retains repetitions. In particular, averaging over it gives equal weight to $G_n$ and $\theta_{t,n}(G_n)$, regardless of their intersection.

Choose pairs $c_1,c_2\in\cU(\mathbb M_n(\C))$ and $d_1,d_2\in\cU(\mathbb M_m(\C))$ as in Proposition~\ref{gap}. Set
\[
N=\mathbb M_n(\C)\otimes(M_n\otimes M_n)\otimes\mathbb M_m(\C).
\]
Let $(e_{ij})$ and $(f_{ij})$ be matrix units for the first and last factors, respectively, and define
\begin{align}
 a_1&=c_1\otimes1\otimes1,
 &a_2&=c_2\otimes1\otimes1,
 &a_3&=\sum_{i=1}^ne_{ii}\otimes(X_{n,i}\otimes1)\otimes1,
 \label{a-triple}\\
 b_1&=1\otimes1\otimes d_1,
 &b_2&=1\otimes1\otimes d_2,
 &b_3&=\sum_{j=1}^m1\otimes V_j\otimes f_{jj}.
 \label{b-triple}
\end{align}
%The use of matrix units here is unrelated to the tensor-length projection $e_l$ above. 
All six elements are unitaries. Their cross commutators vanish except possibly for
\[
[a_3,b_3]=\sum_{i=1}^n\sum_{j=1}^m
e_{ii}\otimes[X_{n,i}\otimes1,V_j]\otimes f_{jj}.
\]
Since this is block-diagonal, Lemma~\ref{almost-commuting} gives
\begin{equation}\label{cross-bound}
\|[a_\alpha,b_\beta]\|\leq4|t|,
\qquad 1\leq\alpha,\beta\leq3.
\end{equation}

For the rest of the paper, take
\begin{equation}\label{parameters}
n_k=k^4,\qquad t_k=\frac1k,\qquad l_k=10k^2,
\qquad k\geq4.
\end{equation}
Write $m_k=2|G_{n_k}|$, and let $N_k$, $a_{k,i}$ and $b_{k,j}$ denote the matrix algebra and unitaries above with $n=n_k$ and $t=t_k$. Notice that $l_k<n_k$ and $l_k/n_k\rightarrow0$.

Let $\alpha_k,\beta_k:\F_3\rightarrow\cU(N_k)$ be the homomorphisms determined by the two triples, and define
\begin{equation}\label{asymptotic}
\pi_k:\F_3\times\F_3\rightarrow\cU(N_k),\qquad
\pi_k(g,h)=\alpha_k(g)\beta_k(h).
\end{equation}
If $g',h$ have reduced word lengths $|g'|,|h|$, then \eqref{cross-bound} gives
\[
\|[\alpha_k(g'),\beta_k(h)]\|
\leq4|g'||h|t_k.
\]
%The same bound holds for generators replaced by their inverses, by unitary invariance of the norm. 
Thus
\[
\|\pi_k(g,h)\pi_k(g',h')-\pi_k(gg',hh')\|
\leq4|g'||h|t_k\longrightarrow0.
\]
In particular, $(\pi_k)_k$ is an asymptotic homomorphism in the normalized $2$-norm.

\section{Proof of Theorem~\ref{main}}\label{proof-section}

Let $u_1,u_2,u_3$ and $v_1,v_2,v_3$ denote the canonical generators of $C^*(\F_3\times\F_3)$. We first show that the asymptotic homomorphism constructed in Section~\ref{construction-section} does not admit u.c.p. approximations.

\begin{proposition}\label{obstruction}
Let $a_{k,i},b_{k,j}\in N_k$, $1\leq i,j\leq3$, be as above. There is no sequence of u.c.p. maps $\varphi_k:C^*(\F_3\times\F_3)\rightarrow N_k$ such that
\begin{equation}\label{assumed-approximation}
\max_{1\leq i,j\leq3}
\{\|\varphi_k(u_i)-a_{k,i}\|_2,
\|\varphi_k(v_j)-b_{k,j}\|_2\}\longrightarrow0.
\end{equation}
\end{proposition}

\begin{proof}
Suppose, toward a contradiction, that such maps exist. Apply Lemma~\ref{state} with $\eps=\eps_k$ equal to the maximum in \eqref{assumed-approximation}. We obtain Hilbert spaces $\cK_k$, unit vectors $\xi_k\in\cK_k\otimes L^2(N_k)$, and positive operators $h_k$ on $L^2(N_k)$ such that
\begin{gather}
\|h_k\|_{\HS}=1,\qquad
\|(1\otimes\Ad(b_{k,j}))\xi_k-\xi_k\|\longrightarrow0,
\label{xi-invariance}\\
\|[h_k,R_{a_{k,i}}]\|_{\HS}\longrightarrow0,
\label{h-invariance}\\
\|(1-e)h_k\|_{\HS}=\|(1\otimes(1-e))\xi_k\|
\quad\text{for every orthogonal projection }e\text{ on }L^2(N_k).
\label{h-xi-support}
\end{gather}
All limits below are ordinary limits as $k\rightarrow\infty$. The notation $o(1)$ denotes a quantity tending to zero. Where no confusion can arise, we suppress $k$ and write $n=n_k$, $t=t_k$ and $m=m_k$.

\smallskip
\noindent\emph{Step 1. Localization using the second free factor.}
The elements $b_1,b_2$ act only on the last factor $\mathbb M_m(\C)$. By \eqref{xi-invariance} and \eqref{gap-vectors}, the distance from $\xi$ to
\[
\cK\otimes L^2(\mathbb M_n(\C))\otimes L^2(M_n\otimes M_n)\otimes\C1_m
\]
is $o(1)$. Let $\zeta$ be the orthogonal projection of $\xi$ onto this space and write $\zeta=\zeta_0\otimes1_m$. Then $\|\zeta-\xi\|=o(1)$. Since $\Ad(b_3)$ is unitary, \eqref{xi-invariance} also gives $\|(1\otimes\Ad(b_3))\zeta-\zeta\|=o(1)$.
%\[ \|(1\otimes\Ad(b_3))\zeta-\zeta\| \leq2\|\zeta-\xi\|+\|(1\otimes\Ad(b_3))\xi-\xi\|=o(1).\]
Writing $1_m=\sum_jf_{jj}$ and using the orthogonality of the $f_{jj}$ in the normalized $2$-norm, we get the identity
\begin{equation}\label{controlled-average}
\|(1\otimes\Ad(b_3))\zeta-\zeta\|^2
=\frac1m\sum_{j=1}^m
\|(1\otimes1\otimes\Ad(V_j))\zeta_0-\zeta_0\|^2.
\end{equation}
The first two identity operators on the right act on $\cK$ and $L^2(\mathbb M_n(\C))$, respectively.

Set $\cL=\cK\otimes L^2(\mathbb M_n(\C))$. The two halves of the list $(V_j)_{j=1}^m$ enumerate $G_n$ and $\theta_{t,n}(G_n)$. By \eqref{group-commutant}, their commutants in $M_n\otimes M_n$ are $A_n\otimes1$ and $\theta_{t,n}(A_n\otimes1)$, respectively. Lemma~\ref{averaging} and \eqref{controlled-average} therefore imply
\begin{align}
\dist\big(\zeta_0,\cL\otimes L^2(A_n\otimes1)\big)&=o(1),
\label{first-localization}\\
\dist\big(\zeta_0,\cL\otimes L^2(\theta_{t,n}(A_n\otimes1))\big)&=o(1).
\label{second-localization}
\end{align}
%Indeed, each half-average in \eqref{controlled-average} is nonnegative, and each is twice the corresponding squared distance.

\smallskip
\noindent\emph{Step 2. Concentration on tensors of low degree.}
Let $\eta$ be the orthogonal projection of $\zeta_0$ onto $\cL\otimes L^2(A_n\otimes1)$ and write
\[
\eta=\sum_{S\subset\{1,\ldots,n\}}\eta_S\otimes(w_S\otimes1),
\qquad \eta_S\in\cL.
\]
By \eqref{first-localization}, $\|\eta-\zeta_0\|=o(1)$. Combining this with \eqref{second-localization} and then applying \eqref{deformation-distance}, we get
\begin{equation}\label{weighted-tail}
\sum_S(1-\rho_t^{2|S|})\|\eta_S\|^2=o(1).
\end{equation}
For $|t|\leq1$, the elementary inequality $\cos t\leq1-t^2/3$ gives $\rho_t\leq1-t^2/6$. Hence, for the parameters in \eqref{parameters},
\begin{equation}\label{tail-gap}
\rho_{t_k}^{2(l_k+1)}
\leq\exp\left(-\frac{l_k+1}{3k^2}\right)
\leq e^{-10/3}<\frac12.
\end{equation}
For $|S|>l_k$, the coefficient $1-\rho_{t_k}^{2|S|}$ in \eqref{weighted-tail} is therefore at least $1/2$. It follows that
\[
\sum_{|S|>l_k}\|\eta_S\|^2=o(1).
\]

Let $e_k$ be the orthogonal projection of $L^2(N_k)$ onto
\begin{equation}\label{low-subspace}
L^2(\mathbb M_{n_k}(\C))\otimes
\Span\{w_S\otimes1\mid |S|\leq l_k\}\otimes\C1_{m_k}.
\end{equation}
The preceding estimate, together with $\|\eta-\zeta_0\|=o(1)$ and $\|\zeta-\xi\|=o(1)$, gives
\begin{equation}\label{xi-low}
\|(1\otimes(1-e_k))\xi_k\|=o(1).
\end{equation}
Using \eqref{h-xi-support} and the self-adjointness of $h_k$, we conclude that
\begin{equation}\label{h-low}
\|h_k-e_kh_ke_k\|_{\HS}
\leq\|(1-e_k)h_k\|_{\HS}+\|e_kh_k(1-e_k)\|_{\HS}
\leq2\|(1-e_k)h_k\|_{\HS}=o(1).
\end{equation}

\smallskip
\noindent\emph{Step 3. Averaging the right action of the first matrix factor.}
Put $H_k=e_kh_ke_k$. The projection $e_k$ commutes with right multiplication by the first matrix factor $\mathbb M_{n_k}(\C)$. In particular, for $i=1,2$,
\[
[H_k,R_{a_{k,i}}]=e_k[h_k,R_{a_{k,i}}]e_k,
\]
so these commutators tend to zero in Hilbert-Schmidt norm by \eqref{h-invariance}.

Let $\mathcal E_k$ be the map on $\HS(L^2(N_k))$ obtained by averaging conjugation by the right action of $\cU(\mathbb M_{n_k}(\C))$ on the first factor, and put $Y_k=\mathcal E_k(H_k)$. The opposite-algebra version of \eqref{gap-operators} gives
\[
\|H_k-Y_k\|_{\HS}^2
\leq\kappa^2\sum_{i=1}^2\|[H_k,R_{a_{k,i}}]\|_{\HS}^2=o(1).
\]
Together with \eqref{h-low}, this yields
\begin{equation}\label{Y-properties}
\|Y_k-h_k\|_{\HS}=o(1),\qquad
\|Y_k\|_{\HS}\longrightarrow1,\qquad
\|[Y_k,R_{a_{k,3}}]\|_{\HS}\longrightarrow0.
\end{equation}
%For the last assertion, use $\|[Y_k,R_{a_{k,3}}]\|_{\HS}\leq2\|Y_k-h_k\|_{\HS}+\|[h_k,R_{a_{k,3}}]\|_{\HS}$.
Moreover, since $e_k$ commutes with the whole right action being averaged,
\begin{equation}\label{Y-support}
Y_k=e_kY_ke_k.
\end{equation}

\smallskip
\noindent\emph{Step 4. The Fourier contradiction.}
Consider the larger subspace
\[
\mathscr H_k^0=L^2(\mathbb M_{n_k}(\C))\otimes
L^2(A_{n_k}\otimes1)\otimes\C1_{m_k}\subset L^2(N_k).
\]
It reduces $R_{a_{k,3}}$: multiplication by each $X_{n_k,i}\otimes1$ preserves $A_{n_k}\otimes1$, and $a_{k,3}$ is self-adjoint. By \eqref{Y-support}, $Y_k$ vanishes on $(\mathscr H_k^0)^\perp$ and has range in $\mathscr H_k^0$. We may therefore compute the commutator in \eqref{Y-properties} on $\mathscr H_k^0$.

Under the unitary identification
$L^2(\mathbb M_{n_k}(\C))=\C^{n_k}\otimes\overline{\C^{n_k}}$
given by $\sqrt{n_k}e_{ij}\mapsto e_i\otimes\overline{e_j}$,
right multiplication acts on the second tensor factor. Reordering the factors identifies
\[
\mathscr H_k^0
=\overline{\C^{n_k}}\otimes
\big(\C^{n_k}\otimes L^2(A_{n_k})\big).
\]
Since $Y_k$ commutes with the full right matrix algebra, there is an operator $y_k$ on $\C^{n_k}\otimes L^2(A_{n_k})$ such that
\begin{equation}\label{Y-factorization}
Y_k|_{\mathscr H_k^0}=\frac{I_{n_k}}{\sqrt{n_k}}\otimes y_k,
\qquad \|y_k\|_{\HS}=\|Y_k\|_{\HS}\longrightarrow1.
\end{equation}
The normalization follows from $\|I_{n_k}\|_{\HS}=\sqrt{n_k}$. Furthermore, \eqref{Y-support} says exactly that
\[
y_k=(1\otimes f_{l_k})y_k(1\otimes f_{l_k}).
\]
On this decomposition, right multiplication by $a_{k,3}$ is
$R_{a_{k,3}}|_{\mathscr H_k^0}
=\sum_{i=1}^{n_k}e_{ii}\otimes(1\otimes M_i)$,
where $e_{ii}$ is the $i$-th coordinate projection on $\overline{\C^{n_k}}$ and $M_i$ is multiplication by $X_{n_k,i}$ on $L^2(A_{n_k})$. Orthogonality of these matrix blocks gives
\begin{equation}\label{controlled-fourier}
\|[Y_k,R_{a_{k,3}}]\|_{\HS}^2
=\frac1{n_k}\sum_{i=1}^{n_k}\|[y_k,1\otimes M_i]\|_{\HS}^2.
\end{equation}
By Lemma~\ref{fourier}, the right side is at least
\[
2\left(1-2\sqrt{\frac{l_k}{n_k}}\right)\|y_k\|_{\HS}^2
\longrightarrow2,
\]
since $l_k/n_k=10/k^2\rightarrow0$. This contradicts \eqref{Y-properties} and proves the proposition.
\end{proof}

\begin{proof}[Proof of Theorem~\ref{main}]
Assume first that $C^*(\F_3\times\F_3)$ has the LLP. The maps $\pi_k$ from \eqref{asymptotic} form an asymptotic homomorphism into the matrix algebras $N_k$. Applying Corollary~\ref{ISW}, we obtain u.c.p. maps $\varphi_k:C^*(\F_3\times\F_3)\rightarrow N_k$ such that
\[
\|\varphi_k(u_i)-a_{k,i}\|_2\longrightarrow0,
\qquad
\|\varphi_k(v_j)-b_{k,j}\|_2\longrightarrow0,
\qquad 1\leq i,j\leq3.
\]
This contradicts Proposition~\ref{obstruction}. Thus $C^*(\F_3\times\F_3)$ does not have the LLP.
Since $\F_3\times\F_3$ is a subgroup of $\F_2\times\F_2$, the subgroup permanence of the LLP \cite[Remark~1.2]{ISW20} implies that $C^*(\F_2\times\F_2)$ cannot have the LLP either. This proves Theorem~\ref{main}.
\end{proof}

\begin{proof}[Proof of Corollary~\ref{products}]
The assumptions imply that $\F_2\times\F_2$ embeds into $\Gamma_1\times\Gamma_2$. Apply Theorem~\ref{main} and the subgroup permanence of the LLP \cite[Remark~1.2]{ISW20}.
\end{proof}

%\clearpage
\small{
}

\end{document}